\documentclass{amsart}

\numberwithin{equation}{section}

\usepackage{amsmath}
\usepackage{amssymb}
\usepackage{enumerate, xspace}
\usepackage{hyperref}

\newtheorem{thmm}{Theorem}
\newtheorem{thm}{Theorem}[section]
\newtheorem{lem}[thm]{Lemma}
 
\newtheorem{sublem}[thm]{Sub-lemma}
\newtheorem{prop}[thm]{Proposition}

\newtheorem{defin}{Definition}

\newtheorem{rem}[thm]{Remark}

\newcommand\B{{\mathcal B}}
\newcommand\Co{{\mathcal C}}
\newcommand\D{{\mathcal D}}

\newcommand\Lp{{\mathcal L}}
\newcommand\Or{{\mathcal O}}

\newcommand\T{{\mathcal T}}

\newcommand\BB{{\mathbb B}}
\newcommand\C{{\mathbb C}}

\newcommand\N{{\mathbb N}}
\newcommand\R{{\mathbb R}}
\newcommand\To{{\mathbb T}}
\newcommand\Z{{\mathbb Z}}

\newcommand\ve{\varepsilon}
\newcommand\vf{\varphi}

\newcommand\de{d^\flat}
\newcommand\Id{\text{\bf Id}}
\newcommand\Tr{\text{\bf Tr\,}}
\newcommand\Trs{\operatorname{Tr}}
\newcommand\Exp{\operatorname{Exp}}

\newcommand{\coef}{\Gamma}
\newcommand{\Diff}{\operatorname{Diff}}

\begin{document}

\title[Fredholm determinants]{Fredholm determinants, Anosov maps and Ruelle resonances}

\author[Carlangelo Liverani]{ }
\address{Carlangelo Liverani\\
Dipartimento di Matematica\\
II Universit\`{a} di Roma (Tor Vergata)\\
Via della Ricerca Scientifica, 00133 Roma, Italy.}
\email{{\tt liverani@mat.uniroma2.it}}


\thanks{It is a real pleasure to thank V.Baladi, D.Dolgopyat,
S.Gou\"ezel for very helpful discussions. In addition, I like to thank
an anonymous referee for suggesting me to address Proposition \ref{prop:ruelle}.} 
\date{April 25, 2005}
\begin{abstract}
{I show that the dynamical determinant, associated to an
Anosov diffeomorphism, is the Fredholm 
determinant of the corresponding Ruelle-Perron-Frobenius transfer
operator acting on appropriate Banach spaces. 
As a consequence it follows, for example, that the zeroes of the dynamical determinant
describe the eigenvalues of the transfer operator and the Ruelle
resonances and that, for $\Co^\infty$ Anosov diffeomorphisms, the
dynamical determinant is an entire function.}
\end{abstract}
\maketitle

\section{Introduction}
In the last years there has been a considerable interest in the study
of dynamical determinants and dynamical
zeta functions (see \cite{BaBa, Go, Ba1, BJR, Ru3, BuK, Na,  Ki, BaRu,
Po1}, just to mention a 
few, \cite{Ru3, Po2} for brief reviews of the field,
\cite{Babook} for a general introduction and \cite{chaos} for a
detailed discussion of physics related issues).  Here, I will focus
on Anosov diffeomorphisms $T$ and the associated Fredholm
determinant $\de$ \eqref{eq:zetaf} for the transfer operator
$\Lp$ (see \eqref{eq:transfer} for a precise definition).

The most satisfactory results have been obtained for analytic
systems \cite{Ru1, Fr1, Fr3, Rug2} and $\Co^{r+1}$ 
expanding maps \cite{Ru2, Fr2}. For axiom A analytic and
$\Co^\infty$-expanding maps the
above mentioned papers prove that the dynamical determinant is an
entire function and its zeroes are exactly the inverse of the
eigenvalues of the associated transfer operator; that is, it can be
interpreted as a Fredholm determinant.

On the contrary for $\Co^{r+1}$ Axiom A  maps or flows the situation is
still unsatisfactory. The strongest
result to date is \cite{Ki} where it is showed that the dynamical determinant
for a $C^{r+1}$ Anosov map, with expansion and contraction estimated by
$\lambda$, is analytic in the disk $\lambda^{\frac r2}$. Nevertheless, in
\cite{Ki} the relation between the dynamical determinant and the
transfer operator is only a formal one, in 
particular no information is available concerning the relation
between the zeroes of such a function and the spectrum of the transfer
operator. It was therefore a bit arbitrary to call such a function a
Fredholm determinant. 

In the present paper the missing relation is derived at the price
of establishing the result in a smaller disk. Building on the results
in \cite{GL} I  will show that it is possible to make sense of the
na\"\i ve idea of smoothing the singular kernel of the transfer
operator, \cite[page 103]{Babook}. This yields a strategy greatly
simplified with respect to previous approaches. In fact, it
essentially boils down to a couple of pages computation. As a
consequence one establishes the complete description of the
correlation spectra (Ruelle resonances) in terms of periodic
orbits. Finally, let me remark that, most likely,  the present approach can be
extended to more general transfer operators (e.g., with smooth weights), systems (e.g., Axiom A) and
to the study of dynamical zeta functions (since the latter can be
expressed as ratios of dynamical determinants \cite{Ru1}).

The plan of the paper is as follows. Section \ref{sec:one} details
and proves the main results of the paper. Given the existence of a
scale of {\em adapted Banach spaces} (see Definition
\ref{def:adapted}) and Lemma \ref{lem:key} the proofs 
are completely self-contained. Lemma \ref{lem:key} is proven in section
\ref{sec:proofs} while Proposition \ref{prop:basic}, proven in
section \ref{sec:cond}, states the existence of the adapted spaces. 
This last result relies on a scale of Banach
spaces introduced in \cite{GL}, yet it should be emphasized that other
choices of adapted spaces are possible, e.g. V.Baladi has recently introduced a different choice that, in
very special cases, enjoys some useful extra properties \cite{Ba2} and V.Baladi with M.Tsujii have announced a different possibility that could yield sharper bounds. 
Finally, an appendix contains an hardly surprising technical result
that, for lack of references, needed to be proven somewhere.

\begin{rem}
In this paper $C$ stands for a generic constant depending only on the
dynamical system $(X,T)$ under consideration. Its actual value can
thus change from one occurrence to the next.
\end{rem}

\section{The results}\label{sec:one}
In the following, I will discuss only the case $X=\To^d$ with the
Euclidean metric. This
simplifies the presentation and the notations since one can avoid the need to
introduce local charts. The general case can be investigated in
complete analogy  by using partitions of unity and local charts along
the lines exploited in \cite{GL}. Also, I will discuss only the
transfer operator associated to the SRB measure although I do not see
any real obstacle in treating more general, smooth,
potentials.\footnote{The main problem is than one needs the extension
of \cite{GL} to such a setting. This is rather straightforward but to
include it here would substantially increase the length of the paper
without adding much to the presentation of the basic idea.}

Let $T\in\operatorname{Diff}^{r+1}(X,X)$ and $\D_r'$ be the space of distributions of order $r$, the transfer operator $\Lp:\D_r'\to\D_r'$ is defined by\footnote{Usually, the transfer operator is defined as acting on function but, in the present contest, it turns out to be essential that the operator can be defined also on distributions. In fact, by using the standard identification between functions and distributions (see also footnote \ref{foot:abuse}), on can restrict the operator to $\Co^r$ obtaining,  for each $h\in \Co^r(X,\R)$, the usual formula $\Lp
h(x):=f\circ T^{-1}|\det(D_xT^{-1})|$ which describes the evolution, under the dynamics,
of the density of the measures absolutely continuous with respect to
the Lebesgue measure $m$.} 

\begin{equation}
\label{eq:transfer}
(\Lp h,\phi):=(h,\phi\circ T), \quad \forall \phi\in\Co^r(X,\C).
\end{equation}
In addition, consider a convolution operator $\hat Q_\ve:\Co^\infty\to\Co^\infty$:\footnote{For
a general manifold $X$, one must introduce coordinates charts $\Psi_i$
and a subordinate partition of unity $\{\phi_i\}$, then, for $\ve$
small enough, one can define
\[
\hat Q_\ve f(x):=\sum_i\int
J\Psi_i(\Psi_i^{-1}(x)+\xi)q_\ve\circ\Psi_i(\Psi_i^{-1}(x)+\xi)
\phi_i\circ\Psi_i(\Psi_i^{-1}(x)+\xi)f\circ\Psi_i(\Psi_i^{-1}(x)+\xi) 
\]
and the following holds essentially unchanged.}
\begin{equation}
\label{eq:convol}
\begin{split}
&\hat Q_\ve f(x):=\int_{\R^d} q_\ve(x-y)f(y) dy\\
&\int_{\R^d} q_\ve(x)=1;\quad \int_{\R^d} x^\alpha q_\ve(x)=0.
\end{split}
\end{equation}
for each multi-index $\alpha$ such that $0<|\alpha|\leq r$, and where
$q_\ve(x)=\ve^{-d}\bar q(\ve^{-1}x)$, $\bar q(x)=\bar q(-x)$, $\text{supp}\,\bar q\subset
\{x\in\R^d\;:\;|x|\leq 1\}$, $\bar q\in\Co^\infty$. 

By duality one can then define the operator $Q_\ve:=\hat Q_\ve':\D'\to\D'$ which can be easily seen to be an extension of $\hat Q_\ve$ to the space of distributions.

It is well known that the spectral properties of the transfer operator depend drastically on the space on which it acts. The space of distributions turns out to be too large of a space to be useful, yet it is well known that $\Co^\infty$ is by far too small and the spectra of $\Lp$ on such a space bears little relevance on the statistical properties of the system. Below we give an abstract characterization of some properties that \emph{ good} dynamical spaces should enjoy.

\begin{defin}
\label{def:adapted}
Given $T\in\operatorname{Diff}^{r+1}(X,X)$, a scale
of Banach spaces $\{\B^s\}_{s\in\N}$ is called {\em adapted to $T$}
if there exists $s_r\in\N$ such that, for $s\in\{1,\dots,s_r\}$,
$\D'_r\supset\B^{s-1}\supset \B^s\supset
\Co^s$.\footnote{\label{foot:abuse}
Of course, to make sense of such a scale it is necessary to
slightly abuse notations and identify each functions $f\in\Co^s(X,\R)$
with a linear functional (distribution) via the standard duality relation $(f,\vf):=\int_X f\vf dm$.}
More precisely, $\|\cdot\|_{s-1}\leq\|\cdot\|_{s}\leq
|\cdot|_{\Co^s}$, and $\overline{\Co^s}^{\,\|\cdot\|_s}=\B^s$.
In addition, there exist $\theta\in(0,1)$ and
$\ve_1>0$ such that, for all $0<s\leq s_r$, $\Lp\in L(\B^s,\B^s)$ and,
for each $\ve\in (0,\ve_1)$, $l< s$ and $n\in\N$,\footnote{By $D_s$ I mean any derivative in the stable direction.}
\begin{align}
&\left|\int f\phi\right|\leq \|f\|_{\B^0} |D_s^r\phi|_\infty\text{ for each } f,\phi\in\Co^\infty
\label{eq:hypo00}\\
& \Lp:\B^{s}\to \B^{s-1} \text{ is compact}
\label{eq:hypo0}\\ 
&\|\Lp^n h\|_{\B^{s}}\leq B \|h\|_{\B^{s}} \text{ for each } h\in\B^s
\label{eq:hypo1}\\ 
&\|\Lp^n h\|_{\B^{s}}\leq
A\theta^{sn}\|h\|_{\B^{s}}+B\|h\|_{\B^{s-1}} \text{ for each }
h\in\B^s\label{eq:hypo2}\\ 
&Q_\ve\in L(\B^s,\Co^\infty)\text{ and }\|Q_\ve-\Id\|_{\B^{s}\to
\B^{s-l}}\leq D\ve^{l}\label{eq:hypo4}\\
&\|hf\|_{\B^s}\leq D'\|h\|_{\B^s}|f|_{\Co^{r}}  \text{ for each }
h\in\B^s,\, f\in\Co^{r}\label{eq:hypo5}
\end{align}
where $A, B, D, D'$ do not depend on $\ve$ and $n$.\footnote{In
fact, Property \eqref{eq:hypo5}, is needed only in the proof of Proposition
\ref{prop:ruelle}.}
If $T\in\Diff^\infty(X,X)$ and we have an adapted scale for each $r$,
with $\lim_{r\to\infty}s_r=\infty$, then we say that we have a {\em
complete series} of adapted Banach spaces.
\end{defin}
From simple arguments (see, e.g., \cite{GL}) follows that on such spaces the spectrum of $\Lp$ has a
{\em physical} interpretation: it describes the rate of decay of
correlations and it is stable with respect to a large family of
perturbations. In addition, in section \ref{sec:cond} I prove:
\begin{prop}\label{prop:basic}
If $T\in\operatorname{Diff}^{r+1}(X,X)$ is Anosov,\footnote{Anosov
means that there exists a continuous splitting $E^u\otimes E^s$,
$\operatorname{dim}(E^s)=d_s$ and $\operatorname{dim}(E^u)=d_u$, of the
tangent bundle and a constant $\lambda>1$ such that $\|D_xT^n
v\|\leq \lambda^{-n}\|v\|$ for all $v\in E^s(x)$, $x\in X$ and $\|D_xT^{-n}
v\|\leq \lambda^{-n}\|v\|$ for all $v\in E^u(x)$, $x\in X$.} then there
exists a scale of Banach spaces adapted to $T$
with $s_r=\lceil\frac {r-1}{2}\rceil$.\footnote{Given $a\in\R$,
$\lceil a\rceil$ stands for the largest integer $n\leq a$.} If
$T\in\Diff^\infty(X,X)$ then the latter constitute a complete series
of adapted spaces. 
\end{prop}

\begin{rem}\label{rem:hypo}
In section \ref{sec:cond} I will define the wanted spaces based on the Banach spaces introduced in \cite{GL}, yet the present results hold for any other choice of adapted Banach spaces satisfying \eqref{eq:hypo00}-\eqref{eq:hypo5}.  
\end{rem} 

\begin{rem}\label{rem:spect} 
The existence of an adapted scale of Banach spaces implies (see
\cite{Babook}):  For each $1<s\leq s_r$ and  $\sigma\in
(\theta,1)$, the operator $\Lp$ is quasicompact on $\B^s$, more
precisely  it can 
be decomposed as $\Lp=P_{\sigma,s}+R_{\sigma,s}$  
where $P_{\sigma,s} R_{\sigma,s}=R_{\sigma,s} P_{\sigma,s}=0$,
$P_{\sigma,s}$ is of finite rank and
\begin{equation}
\label{eq:rest}
\|R_{\sigma,s}^n\|_{\B^s}\leq C\sigma^{sn}.
\end{equation}
\end{rem}

For further use let us set 
\begin{equation}
\label{eq:coeff}
\coef_n:=\sum_{x\in\text{Fix }T^n} |\det(\Id-D_{x}T^n)|^{-1},
\end{equation}

The following estimate is more or less standard. The proof can be
found at the end of the section and is enclosed only for completeness.
\begin{lem}\label{lem:zeroest}
If $T\in\operatorname{Diff}^{r+1}(X,X)$ is Anosov,  for all $n\in\N$,
$\coef_n\leq C$.  
\end{lem}

The main result of the paper is the following.
\begin{thmm}\label{thm:main}
If $T\in\operatorname{Diff}^{r+1}(X,X)$ is Anosov, define, for $|z|<1$,
\begin{equation}\label{eq:zetaf}
\de(z):=\Exp\left(-\sum_{n=1}^\infty\frac{z^{n}}n
\sum_{x\in\text{Fix }T^n} |\det(\Id-D_{x}T^n)|^{-1}\right),
\end{equation}
and consider its analytic extension.
Then, if $\bar s_r:=\lceil\frac{s_r}{4(1+d/s_r)}\rceil>0$,
$\de(z)\det(\Id-zP_{\sigma,\bar s_r})^{-1}$  
is holomorphic and never zero in the disk
$|z|<\theta^{-\bar s_r}$.
Thus in such a disk $\de$ is holomorphic and its zeroes are in one one correspondence with
the eigenvalues of the operator $\Lp$. In addition, the algebraic multiplicity of the zeroes equal the dimension of the associated eigenspaces. In particular, if
$T\in\Co^\infty$, then $\de$ is an entire function. 
\end{thmm}
\begin{rem}
Using the spaces in \cite{GL}, see section \ref{sec:cond}, 
the above Theorem gives the analyticity domain
$\lambda^{\frac{r-1}{8(1+d/(r-1))}}$. This is certainly far from  optimal
and can be easily improved.\footnote{For example, one
should be able to prove $Q_\ve\in L(\B^s,\B^s)$, $d$ can be replaced
by $d_u$, etc.} Yet, since I do not see how to obtain the
expected Kitaev-like bound $\lambda^{-\frac r2}$, I will not strive
for superficial 
improvements at the expenses of clarity, simplicity and brevity.
\end{rem}

In the present language the SRB measures are the eigenvectors
associated to the eigenvalue one. Probably the most interesting
physical consequence of Theorem \ref{thm:main} is the following.
\begin{prop}\label{prop:ruelle}
Given a transitive (hence mixing) Anosov map $T\in\Diff^\infty(X,X)$, let $\mu_{SRB}$
be the SRB measure and let $f,g\in\Co^\infty(X,\R)$, $\mu_{SRB}(f)=\mu_{SRB}(g) =0$, then the function (the  {\em correlation spectra}) 
\[
C_{f,g}(e^{i\omega})=\sum_{n\in\Z}e^{i\omega n}\mu_{SRB}(f g\circ T^n)\;;
\quad \omega\in\R
\]
extends to a meromorphic function on $\C\setminus\{0\}$ and its poles
(often called {\em Ruelle resonances}) are exactly described by the zeroes of $\de$. 
\end{prop}
\begin{proof}
Fix any $r\in\N$ and consider an associated scale of adapted spaces.
Calling $m$ the Lebesgue measure\footnote{Note that, since $(1,\phi)=\int \phi dm$ (see footnote 
\ref{foot:abuse}), then $m$ can also be seen as the element $1$ of $\B^s$ or $\D'_s$.} the SRB measure can be defined as
\[
\mu_{SRB}(\phi)=\lim_{n\to\infty}m(\phi\circ T^n)=\lim_{n\to\infty}\Lp^n 1(\phi).  
\]
Note that, since the map is mixing, then one must be a simple eigenvalue and no other eigenvalues can be present on the unit circle.\footnote{Indeed, if $e^{i\theta}\in\sigma(\Lp)$, then 
$\Pi_\theta:=\lim\limits_{n\to\infty} \frac 1n \sum_{k=0}^{n-1}e^{-i\theta k}\Lp^k$ is well defined and is exactly the projector on the associated eigenspace. Moreover, from \eqref{eq:hypo1} and \eqref{eq:hypo2} follows that $\text{Range}(\Pi_\theta)\subset\B_0$. Hence, by \eqref{eq:hypo00},
\[
|(\Pi_\theta h,\phi)|=|(\Pi_\theta^2 h,\phi)|\leq \lim_{n\to\infty}\frac 1n\sum_{k=0}^{n-1}|(\Pi_\theta h,\vf\circ T^n)\leq \|\Pi_\theta h\|_s|\phi|_\infty.
\]
That is the eigenspace would consist of measures, whereby violating the mixing assumption.}
Consequently, remembering Remark \ref{rem:spect}, $\Lp$ has a spectral gap, hence there exists $\rho>1$ such that, for each $h\in\B^s$
 \[
 \|\Lp^n h-(h,1)\mu_{SRB}\|_s\leq C\|h\|_s \rho^{-n}.
 \]
 Since, $\mu_{SRB}(f g\circ T^n)=\lim\limits_{p\to\infty} (\Lp^p 1,f g\circ T^n)$ it is natural to define  the measures $m_{p,f}(h):=(\Lp^p 1,fh)$. In fact,  by \eqref{eq:hypo5}, we have
$m_{p,f}\in\B^s$, $s\leq s_r$ and
\[
m_{p,f}(g\circ T^n)=\Lp^nm_{p,f}(g)=\mu_{SRB}(g)m_{p,f}(1)+\Or(\rho^{-n})=\Or(\rho^{-n}).
\]
The above means that, for each $|z|<\rho$, 
\begin{equation}\label{eq:defres}
\begin{split}
  \sum_{n\in\N}z^n\mu_{SRB}(f g\circ T^n)&=\lim_{p\to\infty}\sum_{n\in\N}z^n\Lp^nm_{p,f}(g)\\
&=\lim_{p\to\infty}[(\Id-z\Lp)^{-1}m_{p,f}](g)\\
&=[(\Id-z\Lp)^{-1}\lim_{p\to\infty}m_{p,f}](g)\\
&=[(\Id-z\Lp)^{-1}\mu_f](g)=:G_{f,g}(z),
\end{split}
\end{equation}
where $\mu_f(h):=\mu_{SRB}(fh)$. Since $\mu_f\in\B^s$, for each $s\leq s_r$, by Remark \ref{rem:spect}
follows that the function  
$G_{f,g}$ can be extended to a meromorphic function on
$\{z\in\C\;:\;|z|<\theta^{-s_r}\}$. On the other hand, if $|z|>1$,
\[
 \sum_{n\in\N}z^{-n}\mu_{SRB}(f g\circ T^{-n})=
 \sum_{n\in\N}z^{-n}\mu_{SRB}(g f\circ T^n)=G_{g,f}(z^{-1}).
\]
Hence the formula
\[
C_{f,g}(e^{i\omega})=G_{f,g}(e^{i\omega})+G_{g,f}(e^{-i\omega}) -\mu_{SRB}(fg)
\]
together with \eqref{eq:defres}, shows that $C_{f,g}$ is meromorphic in the annulus
$\{z\in\C\;:\;\theta^{s_r}<|z|<\theta^{-s_r}\}$. By Theorem
\ref{thm:main} its poles are the inverse of the zeroes of $\de$ in the
annulus $\{z\in\C\;:\; \theta^{\bar s_r}<|z|<\theta^{-\bar s_r}\}$. The
Lemma easily follows since we have a complete series of spaces and $r$
can be chosen arbitrarily.
\end{proof}
\begin{rem}
For $\Co^{r+1}$ maps the above argument shows that the correlation function is meromorphic in the anulus $\{\lambda^{-\lceil\frac{r-1}2\rceil}<|z|<\lambda^{-\lceil\frac{r-1}2\rceil}\}$, but the relations between the poles and the zeroes of the dynamical determinat can be established only in the smaller anulus $\{\lambda^{-\lceil\frac{r-1}{8(1+\frac {2d}{r-1})}\rceil}<|z|<\lambda^{-\lceil\frac{r-1}{8(1+\frac {2d}{r-1})}\rceil}\}$ .
\end{rem}
\begin{rem}
Since $C_{f,g}(e^{i\omega})$ is the Fourier
transform of the correlation function, it is a
physically accessible function. Its poles on the complex plane (the
Ruelle resonances) can be computed, e.g. via Pade
approximants, hence they are physically observable as well.
\end{rem}

The proof of Theorem \ref{thm:main} rests on the next basic estimate proven
in section \ref{sec:proofs}.
\begin{lem}\label{lem:key}
For each $n\in\N$, $\sigma\in(\theta,1)$ and $\bar s_r=\lceil
\frac{s_r-1}2\rceil$, holds 
true\footnote{Since $P_{\sigma,s}$ is a finite rank operator, the
usual trace $\Trs$ and the determinant are well defined.}
\[
|\coef_n-\Trs P_{\sigma,\bar s_r}^n|\leq C_\sigma\sigma^{\frac{s_r}{4(1+\frac{d}{s_r})}n}.
\]
\end{lem}

\begin{proof}[Proof of Theorem \ref{thm:main}]
For $|z|<1$, let
\[
g(z):=\det(\Id-zP_{\sigma,s})^{-1}\de(z)=\Exp\left(-\sum_{n=0}^\infty
\frac{z^n}n(\coef_n-\Trs P_{\sigma,s}^n)\right).
\]
By the estimate of Lemma \ref{lem:key} , $g$ is analytic and different
from zero, in the disk $|z|<\sigma^{-\frac{s_r}{4(d+1)}}$. Since
\[
\de(z)=\det(\Id-zP_{\sigma,s})g(z)
\]
the theorem trivially follows from the arbitrariness of $\sigma$.
\end{proof}

\begin{proof}[Proof of Lemma \ref{lem:zeroest}]
Clearly we must worry only about large $n$.
Consider $x\in\text{Fix }T^n$, choose a coordinate system $(\xi,\eta)$
in a neighborhood of $x$ such that $W^u(x)=\{(\xi,0)\}$ and
$W^s(x)=\{(0,\eta)\}$. In such coordinates
\[
D_xT^n=\begin{pmatrix}
        A_n(x)&0\\
        0&B_n(x)
        \end{pmatrix}
\]
where $\|A_n(x)^{-1}\|\leq C\lambda^{-n}$ and $\|B_n(x)\|\leq
C\lambda^{-n}$. Accordingly 
\begin{equation}\label{eq:onest}
\begin{split}
|\det(\Id-D_{x}T^n)| &=|\det A_n(x)\det(\Id-A_n(x)^{-1})
\det(\Id-B_n(x))|\\
&\geq C^{-1}|\det A_n(x)|
= C^{-1}|\det (D_{x}T^n\big|_{E^u})|.
\end{split}
\end{equation}
Let us now consider a small fixed $\rho>0$ and let $W^{u,s}_\rho(z)$
be the unstable and stable manifolds of $z$ of size $\rho$,
respectively. By \eqref{eq:onest} and standard distortion arguments
\[
|\det(\Id-D_{x}T^n)|^{-1}\leq C\rho^{-d_u}\int_{W^u_{\rho}(x)}|\det
(DT^n\big|_{E^u})| d\xi= C\rho^{-d_u}\int_{T^{-n}W^u_{\rho}(x)}d\xi.
\]
Next, let us consider the sets $Z_{\rho}(x):=\cup_{y\in
T^{-n}W^u_\rho(x)} W^s_\rho(y)$ and notice that $\text{Fix }T^n\cap
Z_{2\rho}(x)=\{x\}$. Indeed, let $z\in \text{Fix }T^n\cap
Z_{2\rho}(x)$, then, if $\rho$ has been chosen small enough,
$W^s_{2\rho}(z)\cap W^u_{2\rho}(x)$ contain only one point, let it be
$y$. But, by construction $y\in  W^s_{2\rho}(z)\cap T^{-n}W^u_{2\rho}(x)$,
hence $T^n y\in W^s_{2\rho}(z)\cap W^u_{2\rho}(x)$, that is $T^ny=y$. But
$y=\lim_{n\to\infty}T^{-n}y=x$ and $y=\lim_{n\to\infty}T^{n}y=z$, hence
$x=y=z$.

The above discussion implies that if $x_1,x_2\in\text{Fix }T^n$,
$x_1\neq x_2$, then $Z_\rho(x_1)\cap Z_\rho(x_2)=\emptyset$. Hence
\[
\sum_{x\in\text{Fix }T^n} |\det(\Id-D_{x}T^n)|^{-1}\leq
C\rho^{-d} \sum_{x\in\text{Fix }T^n} m(Z_\rho(x))\leq C\rho^{-d} m(X).
\]
\end{proof}

\section{proof of Lemma \ref{lem:key}}\label{sec:proofs}

The first step in the proof of Lemma \ref{lem:key} is to
define, given an integral operator $Kh(x):=\int \kappa(x,y)h(y)dy$,
$\kappa\in\Co^0(X^2)$,\footnote{Notice that the definition below may
not coincide necessarily with the usual trace even when the latter is
well defined, e.g. it does not necessarily correspond to the
sum of the eigenvalues.}
\begin{equation}\label{eq:trace}
\Tr K:=\int_X \kappa(x,x) dx.
\end{equation}
The first key ingredient is a
representation of such an integral trace for small $\ve$.
\begin{lem}\label{lem:ve0}
For each $\ve$ small enough, holds true
\[
\Tr Q_\ve\Lp^{n}= \coef_n+{\mathcal O}(\ve^{r+1}).
\]
\end{lem}
\begin{proof}
Since
$Q_\ve\Lp^{n}h(x)=\int_{X}q_\ve(x-T^ny) h(y)dy=:\int_X \kappa_{\ve,n}(x,y)h(y)dy$, we have
\begin{equation}
\label{eq:trace1}
\Tr Q_\ve\Lp^{n}=\int_X \kappa_{\ve,n}(x,x)dx=\int_{X}q_\ve(x-T^nx) dx. 
\end{equation}
Next, we consider the change of variable $z=\Phi_n(x):=x-T^nx$,
clearly
\begin{equation}
\label{eq:compdet}
\det D_x\Phi_n=\det (\Id-D_{x}T^n)
\end{equation}

Let $B(0,\ve)\subset\R^d$ be the ball of radius $\ve$ and center zero. If $z\in B(0,\ve)$, then it turns to be useful to define the map $F_z:\To^d\to\To^d$
\begin{equation}
\label{eq:Fz}
F_z(x):=T^n(x)+z\mod 1.
\end{equation}
For $\ve$ small enough, $F_z$ is still hyperbolic hence for each $x\in \text{Fix }(T^n)$ we can consider the $F_z$ $\ve$-orbit $\{x,x,\dots\}$. By shadowing the exists a unique point $x_z$ , in a neighborhood of $x$, such that $x_z=F_z(x_z)=T^n(x_z)+z$. The latter fact means that $\Phi_n(x_z)=z$, that is $B(0,\ve)\subset\text{Range }\Phi_n$. On the other hand,
If  $x\in \Phi^{-1}_n(B(0,\ve))$, then, by shadowing, it is associated to a unique periodic point of
period $n$. Indeed, given $z\in B(0,\ve)$ and $x$ such that $x-T^nx=z$, then
there exists a unique periodic orbit of period $n$ in a neighborhood of the periodic $\ve$-pseudo-orbit $\{x, Tx,\dots, T^{n-1}x\}$. We can then define the function $\Psi: \Phi_n^{-1}(B(0,\ve))\to\text{Fix }T^n$. For each $x\in\text{Fix }T^n$ we can then define the set $\Delta_x:=\Psi^{-1}(x)$. Due to hyperbolicity of $T$ it is easy to verify that $\Phi:\Delta_x\to B(0,\ve)$ is one-to-one beside being onto. We can then label the inverse branches of $\Phi_n$ by the elements of $\text{Fix }T^n$.
\begin{sublem}\label{eq:nastyest}
There exists a constant $M$ such that, for each inverse branch,
\[
\|\Phi_n^{-1}\|_{\Co^{s+1}}\leq
M \;; \quad\|\det(D\Phi_n^{-1})\|_{\Co^{s}}\leq M
|\det(D_{x_*}\Phi_n^{-1})| \quad \forall s\in\{0,\dots,r\}, 
\]
where $x_*\in \text{Fix }T^n$ labels the inverse branch.
\end{sublem}

The above estimate, whose technical but
straightforward proof is  postponed to the appendix, together
with Lemma \ref{lem:zeroest}, \eqref{eq:convol} and \eqref{eq:trace1} yields
\[
\begin{split}
\Tr Q_\ve\Lp^{n}&=\sum_{x\in\text{Fix }T^n}
\int_{\R^d}q_\ve(z) 
\left|\det(D_{\Phi|_{\Delta_x}^{-1}(z)}\Phi_n^{-1})\right|
dz\\
&=\sum_{x\in\text{Fix }T^n} 
|\det(\Id-D_{x}T^n)|^{-1} +{\mathcal O}(\ve^{r}) .
\end{split}
\]
\end{proof}

Hence,\footnote{The following equalities can be easily verified by
direct computation.} 
\begin{equation}
\label{eq:traceone}
\begin{split}
\Tr Q_\ve\Lp^{n}&=\Tr\Lp^{n} Q_\ve=\lim_{\delta_1\to 0}\Tr
Q_{\delta_1}\Lp^{n}Q_\ve\\
&=
\lim_{\delta_1\to 0}\Tr (Q_{\delta_1} P_{\sigma,s}^{n}Q_\ve) +\lim_{\delta_1\to 0}\Tr
(Q_{\delta_1}R_{\sigma,s}^{n} Q_\ve). 
\end{split}
\end{equation}

Setting $\vf_{\ve, y}(x):=q_\ve(x-y)\in\Co^\infty$, for each $A\in
L(\B^s,\B^s)$ and $f\in \Co^s$ holds\footnote{This follows immediately from the fact
that, if $f_n\to f$ in $L^1$, then $Q_\ve f_n\to Q_\ve f$ in $\Co^s$,
hence $\lim_{n\to\infty}AQ_\ve f_n=AQ_\ve f$. One can then approximate
$f$ by piecewise constant function and compute the corresponding
Riemann sums. Taking the limit and since $y\mapsto \vf_{\ve,y}\in\B^s$ is continuous one recovers the integral on the right which is meant in Bochner sense.} 
\begin{equation}
\label{eq:integra}
AQ_\ve f=\int A\vf_{\ve,y} f(y) dy.
\end{equation}
Using \eqref{eq:integra} and remembering \eqref{eq:hypo00},
\eqref{eq:hypo4}, \eqref{eq:rest} yields
\begin{equation}\label{eq:essb}
\begin{split}
\left|\Tr (Q_{\delta_1}R_{\sigma,s}^{n}Q_\ve)\right|&=\left|\int
(Q_{\delta_1}R_{\sigma,s}^{n}\vf_{\ve,y})(y)\right|\leq 
\|Q_{\delta_1}R_{\sigma,s}^{n}\vf_{\ve,y})\|_{\B^{s-1}}\\
&\leq C\|R_{\sigma,s}^{n}\vf_{\ve,y})\|_{\B^{s}}
\leq C\sigma^{sn}\ve^{-s-d}.
\end{split}
\end{equation}

The last step is given by the following perturbation result.
\begin{lem}
\label{lem:pertur}
There exists $\ve_1>0$ such that, for each $\ve\in(0,\ve_1)$ and
$n\in\N$,
\[
|\Trs P_{\sigma,s}^n-\lim_{\delta_1\to 0}\Tr Q_{\delta_1}
P_{\sigma,s}^nQ_\ve|\leq C\ve^{s_r-s}. 
\]
\end{lem}
\begin{proof}
Since $P_{\sigma,s}$ is finite dimensional the usual trace of
$P_{\sigma,s}^n$, $\Trs P_{\sigma,s}^n$, is well defined. More precisely,
$P_{\sigma,s}h=\sum_i w_i\ell_i(h)$, $w_i\in\B^{s_r}$, $\ell_i\in
(\B^s)'$, and $\Trs P_{\sigma,s}=\sum_i\ell_i(w_i)$. Hence, for each $h\in\Co^\infty$, by \eqref{eq:integra} and
\eqref{eq:hypo4}, 
\[
\begin{split}
(Q_{\delta_1} P_{\sigma,s}^nQ_\ve h)(x)&=\int (Q_{\delta_1} P_{\sigma,
s}^n\vf_{\ve,y} )(x)h(y)dy\\
&=\int_X\sum_i (Q_{\delta_1}w_i)(x) \ell_i(P_{\sigma,s}^{n-1}\vf_{\ve,y})h(y) dy.
\end{split}
\]
Thus, by \eqref{eq:integra} and \eqref{eq:hypo4},
\begin{equation*}
\begin{split}
\lim_{\delta_1\to 0}\Tr Q_{\delta_1} P_{\sigma,s}^nQ_\ve&=\lim_{\delta_1\to 0}\sum_i\ell_i(Q_{\delta_1}P_{\sigma,s}^{n-1}Q_\ve w_i)=\sum_i\ell_i(P_{\sigma,s}^{n-1}Q_\ve w_i)\\
&=\Trs P_{\sigma,s}^n +\sum_i\Or(\|P_{\sigma,s}^{n-1}(Q_\ve-\Id) w_i\|_s)
=\Trs P_{\sigma,s}^n+\Or(\ve^{s_r-s}).
\end{split}
\end{equation*}
\end{proof}
Collecting \eqref{eq:traceone}, \eqref{eq:essb} and Lemma
\ref{lem:pertur} yields
\begin{equation}
\label{eq:opla}
\Tr (Q_\ve\Lp^n)= \Trs
(P_{\sigma,s}^n)+\Or(\sigma^{sn}\ve^{-s-d}+\ve^{s_r-s}). 
\end{equation} 
\begin{proof}[Proof of Lemma \ref{lem:key}]
Lemma \ref{lem:ve0} and \eqref{eq:opla} imply
\[
|\coef_n-\Trs P_{\sigma,s}^n|\leq C
\left(\ve^{r}+\ve^{s_r-s}+\sigma^{sn}\ve^{-s-d}\right). 
\]
Finally, choose $\ve=\sigma^{\frac {s}{s_r+d}n}$, $s=\bar s_r=\lceil \frac
{s_r-1} 2\rceil$, and hence the lemma.
\end{proof}

\section{proof of proposition \ref{prop:basic}}\label{sec:cond}

Let us start by recalling the
scale of Banach spaces introduced in \cite{GL}.

It is well known that being Anosov is equivalent to the existence of a
continuous strictly invariant vector field $\Co$. Let $\Co'$ be another
continuous cone field contained in $\text{Int}(\Co)$.  Consider
$\delta>0$ and a set of $d_s$-dimensional manifolds (with boundary)
$\Omega$ such that, if $W\in\Omega$, then there exists $x_W\in W$ and a $d_s$
dimensional hyperplane $E_W$ contained in $\Co'$ such that, making an
isometric change of coordinates such that
$E_W=\{(\xi,0)\;:\;\xi\in\R^{d_s}\}$,
$W=\{x_W+(\xi,\gamma_W(\xi))\;:\;\xi\in\R^{d_s}, \|\xi\|\leq \delta\}$, 
$\T_x W\subset \Co(x)$, for each $x\in W$, and
$|\gamma_W|_{\Co^{r+1}(\R^{d_s},\R^{d_u})}\leq M$ for some fixed $M$
large enough.

Given $W\in \Omega$, we will denote by $\Co_0^q(W,\R)$ the set of functions
from $W$ to $\R$ which are $\lceil q\rceil$ times continuously differentiable and such that the  $\lceil q\rceil$ derivative is $q- \lceil q\rceil$ H\"older continuos
on $W$ and vanish on a neighborhood of the boundary of $W$. For
each $h\in\Co^\infty(X,\R)$ and $q\in\R_+$, $p\in\N^*$, let
  \begin{equation}\label{eq:norm-st}
  \|h\|_{p,q}:=\sup_{|\alpha|\leq p}
  \sup_{W \in \Omega}
  \sup_{\substack{\vf\in\Co_0^{q}(W,\R)\\
  |\vf|_{\Co^{q}}\leq 1}}\int_W  \partial^\alpha h
  \cdot \vf,
  \end{equation}
and define the Banach spaces $\B^{p,q}:=\overline{\Co^\infty}^{\,\|\cdot\|_{p,q}}$.

In \cite{GL} it is proven that,
setting $q_r=\frac {r-1}{2}$ and
$\B^s:=\B^{s,q_r}$, $\Lp\in L(\B^s,\B^s)$ and it satisfies
\eqref{eq:hypo00}, \eqref{eq:hypo0}, 
\eqref{eq:hypo1} and \eqref{eq:hypo2}, provided
$s\leq \lceil\frac{r}{2}\rceil=: s_r$.

\begin{lem}\label{lem:randpert}
For each $l<s\leq s_r$ holds $Q_\ve\in L(\B^s,\Co^\infty)$ and
\[
\|Q_\ve-\Id\|_{\B^s\to\B^{l}}\leq C\ve^{s-l}.
\]
\end{lem}
\begin{proof}
Consider $W\in\Omega$.
Let $W_\rho:=\{x_W+(\xi,(1-\rho)\gamma(\xi))\in\R^d\;:\; \|\xi\|\leq 
\delta\}$. Clearly, provided $\delta$ has been
chosen small enough, for each $\rho\in[0,1]$, $W_\rho\in\Omega$ with
$x_{W_\rho}=x_W$, $E_{W_\rho}=E_W$. Then, for each multi-index
$\alpha$, $|\alpha|\leq l$, 
\vskip-12pt
\[
\begin{split}
&\int_{W} \partial^\alpha h \vf=\int_{\|\xi\|\leq
\delta}\partial^\alpha h(\xi,\gamma(\xi))\vf(\xi,\gamma(\xi))
J_W(\xi)d\xi \\
=&\sum_{i=0}^{s-l-1}\frac 1{i!}\int_{\|\xi\|\leq
\delta}\partial^\alpha
\frac{d^i}{d\rho^i}h(\xi,(1-\rho)\gamma(\xi))\big|_{\rho=\ve} \vf(\xi,\gamma(\xi))
J_W(\xi)d\xi\\
+& \int_0^\ve d\rho_1\cdots\int_0^{\rho_{s-l-1}} d\rho_{s-l}\int_{\|\xi\|\leq
\delta} \partial^\alpha\frac{d^{s-l}}{d\rho^{s-l}}
h(\xi,(1-\rho)\gamma(\xi))\big|_{\rho=\rho_{s-l}} \vf(\xi,\gamma(\xi))
J_W(\xi)d\xi \\
=&\sum_{\substack{i=0\\|\beta|=i}}^{s-l-1}\frac
{(-\ve)^i}{i!}\int_{W_\ve}\gamma^\beta\partial^{\alpha+\beta}h \vf_\ve
+ \sum_{|\beta|=s-l}(-1)^{s-l}\int_0^\ve d\rho_1\cdots\int_0^{\rho_{s-l-1}}
d\rho_{s-l}\\
&\ \ \ \ \times \int_{W_{\rho_{s-l}}}\gamma^\beta\partial^{\alpha+\beta}h
\vf_{\rho_{s-l}} 
\end{split}
\]
\[
=\sum_{\substack{i=0\\|\beta|=i}}^{s-l-1}\frac
{(-\ve)^i}{i!}\int_{W_\ve}\gamma^\beta\partial^{\alpha+\beta}h \vf_\ve
+\Or(\|h\|_{s,q_r}\ve^{s-l})
\]
Thus, in order to estimate the norms, it suffices to consider
$\Omega_\ve:=\{W_\ve\;:\;W\in\Omega\}$, that is manifolds uniformly strictly
inside the cone field. Let $W\in\Omega_\ve$, then
\[
\int_W \partial^\alpha Q_\ve h \vf-\int_W \partial^\alpha h\vf
=\int dz\,q_\ve(z)\int_Wdy\,
[\partial^\alpha h(y+z)-\partial^\alpha h(y)]\vf(y).
\]
If $W=\{\xi, \gamma(\xi)\}$, then $W_z=\{(\xi,\gamma(\xi))+z\}\in\Omega$, provided
$z$ is small enough. Thus, for $|\alpha|=l$, remembering \eqref{eq:convol},
\[
\begin{split}
\int_W \partial^\alpha Q_\ve h \vf-\int_W \partial^\alpha
h\vf&=\sum_{|\beta|=s-l} \int
dzq_\ve(z)\int_0^1 dt_1\cdots\int_0^{t_{s-l-1}}
dt_{s-l}\\
&\ \ \ \times \int_{W_{zt_{s-l}}}dy \partial^{\alpha+\beta}
h(y)z^{\beta}\vf_{zt_{s-l}}(y) \leq\|h\|_{s,q_r}\ve^{s-l}.
\end{split}
\]
\end{proof}

Finally, \eqref{eq:hypo5} follows easily from \eqref{eq:norm-st}. Clearly if
$T\in\Diff^\infty(X,X)$ we have a complete series of adapted spaces.
\appendix
\section{Proof of Sub-Lemma \ref {eq:nastyest}}

Let us choose a periodic point $x_*\in\text{Fix }(T^n)$ and limit our
considerations to the associated inverse branch, that, by a slight
abuse of notation, I will designate simply by $\Phi_n^{-1}$. 
 Again we will use the map $F_z$ introduced in \eqref{eq:Fz}.
Clearly $D_z\Phi_n^{-1}=(\Id-D_{x_z}T^n)^{-1}=(\Id-D_{x_z}F_z)^{-1}$, where
$\Phi_n(x)=z$.

To study the regularities properties at a given point  $z_0$, small enough, it is convenient to
perform an affine change of coordinates $\Lambda(z_0)$ such that
$\Lambda(z_0)(x_{z_0})=0$, $|(D_x\Lambda)|_{\Co^0}
+|(D_x\Lambda)^{-1}|_{\Co^0}\leq C$, and\footnote{Here, and in the
following, given $I\subset \R^q$ and a function $f$ from $I$ to some
Banach algebra $\BB$, by $|\cdot|_{\Co^p}$ we mean the norm
$\sup_{z\in I}\sum_{0\leq|\alpha|\leq
p}\|\partial^\alpha f(z)\|_{\BB}$ so that $\Co^p(I,\BB)$ is itself a Banach algebra.} 
\begin{equation}
\label{eq:forma1}
D_{\Lambda(z_0)(x_{z})}\widetilde F_{z,z_0}=\begin{pmatrix}
        A(z,z_0)&B(z,z_0)\\
        C(z,z_0)&D(z,z_0)
        \end{pmatrix}
\end{equation}
where $\widetilde F_{z,z_0}:=\Lambda(z_0)\circ F_z \circ\Lambda(z_0)^{-1}$;
$\|A(z_0,z_0)\|,\|D(z_0,z_0)^{-1}\|\leq \lambda^{-n}$ and
$B(z_0,z_0)=C(z_0,z_0)=0$. In other words, in the new coordinates, $\{(\xi,0)\}$ corresponds to the stable manifold at $z_0$ and $\{(0,\eta)\}$ to the unstable one.
In such coordinates,\footnote{Here and in the following
I suppress the dependence on $z, z_0$ and $\Lambda$ when non confusion arises.}
\begin{equation}
\label{eq:basicder}
D_z\Phi_n^{-1}\big|_{z=z_0}=
       \begin{pmatrix}
        (\Id-A)^{-1}&0\\
        0&-(\Id-D^{-1})^{-1}
        \end{pmatrix}\begin{pmatrix}
        \Id&0\\
        0&D^{-1}
        \end{pmatrix}.
\end{equation}
Given the simpler structure of  $\tilde F_{z_0,z_0}$ it would be much easier to study its regularity rather than the one of $F_{z_0}$. Yet, the two are equivalent only if the change of coordinates $\Lambda(z_0)$ is uniformly $\Co^r$. To prove the latter is our first task.
 
We start by computing the derivatives of
$x^i_{z}:=T^ix_z$ with respect to $z$: 
\[
\frac{\partial x^i_z}{\partial z}=D_{x_z}T^i(\Id-D_{x_z}T^n)^{-1}.
\]
It is convenient to use in the tangent space at $x^i_{z}$ the
coordinates pushed forward from the tangent space of $x_{z}$. In such
coordinates holds 
\begin{equation}
\label{eq:forma2}
D_{x^{i-1}_{z}}T=:\begin{pmatrix}
        A_i(z,z_0)&B_i(z,z_0)\\
        C_i(z,z_0)&D_i(z,z_0)
        \end{pmatrix}
\end{equation}
where $\|A_i\|,\|D_i^{-1}\|\leq \lambda^{-1}$ and
$B_i(z_0,z_0)=C_i(z_0,z_0)=0$. Hence, 
\begin{equation}
\label{eq:point}
\frac{\partial x^i_z}{\partial z}\big|_{z=z_0}=\begin{pmatrix}
                                  \prod_{j=1}^iA_i &0\\
                                    0 & \prod_{j=i+1}^n D_i^{-1}
                                  \end{pmatrix}
 \begin{pmatrix}
        (\Id-A)^{-1}&0\\
        0&-(\Id-D^{-1})^{-1}
        \end{pmatrix},
\end{equation}
which readily implies $|\partial_zx^i_z|_{\Co^0}\leq \lambda^{-i}C$, for the stable coordinates, and $|\partial_z x^i_z|_{\Co^0}\leq \lambda^{i-n}C$ for the unstable ones. 

An hyperplane $E$ in the stable  direction is uniquely
determined by a linear operator $U:\R^{d_s}\to\R^{d_u}$:
$E=\{(\xi,U\xi)\;:\;\xi\in\R^{d_s}\}$. A simple computation shows
that, by defining
\[
H(z,z_0,U):=(C(z,z_0)+D(z,z_0)U)(A(z,z_0)+B(z,z_0)U)^{-1},
\]
the stable hyperplane for $\widetilde F_z$ at the point $\Lambda(z_0)(x_z)$
is determined by the fixed point of 
$H(z_0,z,U(z))=U(z)$.\footnote{In fact, it is known that $U(z)$  is
$\Co^{r-1}$, e.g. see \cite[Propositions 1, 2]{Po3}, yet here we need explicit
estimates. This forces us to redo the argument.}
 Applying the implicit function theorem, since
by construction $U(z_0)=0$, yields 
\begin{equation}
\label{eq:stableU}
\partial_z U\big|_{z=z_0}=-(\Id-G)^{-1}(D^{-1}\partial_z
C\big|_{z=z_0}), 
\end{equation}
where $G:GL(d_s,d_u)\rightarrow GL(d_s,d_u)$ is defined by
$G(V):=D^{-1}VA$. In addition,
\begin{equation}
\label{eq:Cest}
\partial_{z_i} C\big|_{z=z_0}=\sum_{p=1}^d\sum_{k=1}^n\left[\prod_{j>k}D_j\right]
c_{2,p}(x^{k-1}_{z_0})\left[\prod_{j<k}A_j\right]\frac{\partial
(x^{k-1}_{z})_p}{\partial z_i}\big|_{z=z_0}
\end{equation}
where
\[
\partial_{x_p}D_xT=\begin{pmatrix}
                   a_{2,p}(x)&b_{2,p}(x)\\
                   c_{2,p}(x)&d_{2,p}(x).
                   \end{pmatrix} 
\]
This, by equations \eqref{eq:forma1} and \eqref{eq:forma2} implies
$|U|_{\Co^1}\leq C$. Since the exact same argument can be carried
out for the unstable space, remembering \eqref{eq:point} we have
$|\Lambda|_{\Co^1}+|D_{x_{z_0}}\Lambda|_{\Co^1} +
|D_{x_{z_0}}\Lambda^{-1}|_{\Co^1}\leq C$. We can thus 
conclude the argument by induction: let $l\leq r$ and suppose that
$|x^i_z|_{\Co^l}\leq C$ and $|D_x\Lambda|_{\Co^{l}}
+|D_x\Lambda^{-1}|_{\Co^{l}}\leq C$. Then by 
\eqref{eq:forma2} follows $A_i(z,z), B_i(z,z), C_i(z,z), D_i(z,z)$,
seen as functions of $z$ have $\Co^l$ norms equibounded. In turn, by
\eqref{eq:point}, it follows $|x^i_z|_{\Co^{l+1}}\leq C$. Equations,
\eqref{eq:stableU} and \eqref{eq:Cest}, imply then
$|D_x\Lambda|_{\Co^{l+1}}+|D_x\Lambda^{-1}|_{\Co^{l+1}}\leq C$,
provided $l+1\leq r$. This proves that the change of coordinates $\Lambda(z_0)$ is uniformly $\Co^r$. 

It is now easy to verify that $\tilde F_{z}$ is $\Co^r$ and, remembering \eqref{eq:basicder}, the first inequality of the sub-lemma  can be readily proven. To conclude note that
\[
\begin{split}
\det(D\Phi^{-1}_n)=-\det(\Id-A)^{-1}\det(\Id-D^{-1})\det(D^{-1}).
\end{split}
\]
Next, since given any smooth function $\Delta(z)$ with values in the  invertible matrices,
\[
\begin{split}
&\partial_{z_i}\det \Delta(z)=\lim_{h\to 0}\frac{\det(\Id+[\Delta(z+he_i)-\Delta(z)]\Delta(z)^{-1})-1}{h}\det(\Delta(z))\\
&= \lim_{h\to 0}\frac{e^{\Trs\ln(\Id+[\partial_{z_i}\Delta(z)]\Delta(z)^{-1}h)}-1}h\det(\Delta(z))
=\Trs([\partial_{z_i}\Delta(z)]\Delta(z)^{-1})\cdot\det(\Delta(z)),
\end{split}
\] 
also in view of \eqref{eq:point}, holds true
\[
\|\det(D\Phi^{-1}_n)\|_{\Co^s}\leq C |\det(D^{-1})|.
\]
Finally, since $|\ln\det D|_{\Co^1}\leq
C$, it follows $\|\det(D^{-1})\|_{\Co^{0}}\leq |\det(D^{-1}(x_*))|$
and the lemma.

\end{document}